\documentclass{amsart}
\usepackage{amsmath,amsthm,verbatim,amssymb,amsfonts,amscd,graphicx}
\usepackage{xcolor}
\usepackage{graphics}

\usepackage{euscript, enumerate}
\usepackage{url,hyperref}

\newcommand{\eps}{\varepsilon}

\newcommand{\be}{\begin{equation}}
\newcommand{\ba}{\begin{aligned}}
\newcommand{\bee}{\begin{equation*}}
\newcommand{\ee}{\end{equation}}
\newcommand{\ea}{\end{aligned}}
\newcommand{\eee}{\end{equation*}}
\newcommand{\bea}{\begin{equation} \begin{aligned} }
\newcommand{\eea}{\end{aligned}\end{equation} }

\theoremstyle{plain}
\newtheorem{theorem}{Theorem}[section]

\newtheorem{cor}[theorem]{Corollary}

\theoremstyle{remark}

\theoremstyle{definition}

\numberwithin{equation}{section}

\begin{document}
\title{Uniqueness and stability of singular Ricci flows in higher dimensions}
\author{Robert Haslhofer}

\begin{abstract}
In this short note, we observe that the Bamler-Kleiner proof of uniqueness and stability for 3-dimensional Ricci flow through singularities generalizes to singular Ricci flows in higher dimensions that satisfy an analogous canonical neighborhood property. In particular, this gives a canonical evolution through singularities for manifolds with positive isotropic curvature. The new ingredients we use are the recent classification of higher dimensional $\kappa$-solutions by Brendle, Daskalopoulos, Naff and Sesum, and the maximum principle for the linearized Ricci-DeTurck flow on locally conformally flat manifolds due to Chen and Wu.
\end{abstract}

\maketitle

\section{Introduction}

In a recent breakthrough \cite{BK_uniqueness}, Bamler-Kleiner established uniqueness and stability for 3-dimensional Ricci flow through singularities. This confirmed a long-standing conjecture of Perelman \cite{Per1,Per2}, and had spectacular applications in 3-dimensional geometry and topology \cite{BK_diffeo,BK_contractibility,BK_prime}, specifically the resolution of the generalized Smale conjecture and the contractibility of the space of positive scalar curvature metrics on 3-manifolds.\\

In this short note, we generalize the uniqueness and stability result to singular Ricci flows in higher dimensions that satisfy an analogous canonical neighborhood property. To state our results, recall first that a Ricci flow spacetime $\mathcal{M}=(\mathcal{M},\mathfrak{t},\partial_{\mathfrak{t}},g)$, as introduced by Kleiner-Lott \cite{KL}, is given by an $(n+1)$-dimensional manifold $\mathcal{M}$ with boundary $\partial\mathcal{M}=\mathfrak{t}^{-1}(0)$, a time function $\mathfrak{t}:\mathcal{M}\to [0,\infty)$ without critical points, a time vector field $\partial_{\mathfrak{t}}$ satisfying $\partial_{\mathfrak{t}}\mathfrak{t}=1$, and a metric $g$ on $\ker(d \mathfrak{t})\subset T\mathcal{M}$ satisfying the Ricci flow equation $\mathcal{L}_{\partial_{\mathfrak{t}}}g =-2\mathrm{Rc}(g)$. Now, similarly as in \cite[Definition 3.15]{BK_contractibility} we call a complete Ricci flow spacetime with compact initial condition a singular Ricci flow if for every $\eps>0$ there exists an $r>0$ such that the $\eps$-canonical neighborhood property holds below scale $r$ on $[0,\eps^{-1}]$, i.e. for any $x\in \mathcal{M}\cap\{ \mathfrak{t}\leq \eps^{-1}\}$ with curvature scale less than $r$ the time-slice $(\mathcal{M}_{\mathfrak{t}(x)},g_{\mathfrak{t}(x)},x)$ is $\eps$-close to some time-slice of a $\kappa$-solution. Here, as  in \cite[Definition 1.1]{BN}, a $\kappa$-solution in dimension $n\geq 4$ is an ancient complete nonflat $n$-dimensional Ricci flow that is uniformly PIC and weakly PIC2, has bounded curvature on compact time intervals, and is $\kappa$-noncollapsed at all scales.

\begin{theorem}[uniqueness of singular Ricci flows]\label{thm_uniq}
Any singular Ricci flow is uniquely determined by its initial time-slice $(\mathcal{M}_0,g_0)$ up to isometry. More precisely, if $\mathcal{M}$ and $\mathcal{M}'$ are singular Ricci flows then any isometry $\phi_0:\mathcal{M}_0\to \mathcal{M}_0'$ can be uniquely extended to a smooth diffeomorphism $\phi :\mathcal{M}\to\mathcal{M}'$ such that $\mathfrak{t}=\phi^\ast\mathfrak{t}'$, $\partial_{\mathfrak{t}}=\phi^\ast \partial_{\mathfrak{t}}'$ and $g=\phi^\ast g'$.
\end{theorem}

This generalizes \cite[Theorem 1.3]{BK_uniqueness} to higher dimensions. Concerning applications, recall that Hamilton \cite{Hamilton_pic} and Chen-Zhu \cite{CZ} constructed a Ricci flow with surgery on 4-manifolds with positive isotropic curvature (PIC). By \cite[Appendix A]{KL}, fixing the initial condition, along a sequence of surgery scales $\delta_i\to 0$ one can get as a suitable limit a singular Ricci flow. Together with Theorem \ref{thm_uniq} (uniqueness) we thus obtain:

\begin{cor}[well-posedness of singular Ricci flow with PIC]\label{cor_well}
For any compact 4-manifold $(M_0,g_0)$ with positive isotropic curvature there exists a unique singular Ricci flow with initial condition $(M_0,g_0)$.
\end{cor}

A similar well-posedness result holds for Ricci flow with PIC in dimension $n\geq 12$, where a flow with surgery has been constructed more recently by Brendle \cite{Brendle_pic}. Moreover, it seems likely that the arguments for positive scalar curvature metrics from \cite{BK_contractibility} can be generalized to higher dimensions yielding that the moduli space of PIC metrics on a given compact manifold is contractible (or empty).\\

Similarly as in \cite{BK_uniqueness}, the uniqueness result arises in fact as a consequence of a more quantitative stability result. To state this, recall from \cite[Definition 3.34]{BK_contractibility} that an $\eps$-isometry between Ricci flow spacetimes $\mathcal{M}$ and $\mathcal{M}'$ is a time-preserving diffeomorphism $\phi: \mathcal{U}\to\mathcal{U}'$ between open subsets that are large enough to ensure that the curvature scale is less than $\eps$ on their complements for times less than $\eps^{-1}$, such that the $C^{\lfloor \eps^{-1}\rfloor}$ spacetime norms of $\phi^\ast g'-g$ and $\phi^\ast \partial_{\mathfrak{t}}'-\partial_{\mathfrak{t}}$ are less than $\eps$ for times less than $\eps^{-1}$.

\begin{theorem}[stability of singular Ricci flows]\label{thm_stability}
For every singular Ricci flow $\mathcal{M}$ and every $\eps>0$ there exists a $\delta>0$ with the following significance. If $\mathcal{M}'$ is any singular Ricci flow then any diffeomorphism $\phi_0: \mathcal{M}_0\to \mathcal{M}_0'$ 
satisfying $|| \phi_0^\ast g_0'-g_0||_{C^{\lfloor \delta^{-1}\rfloor}(M)}<\delta$ can be extended to an $\eps$-isometry between $\mathcal{M}$ and $\mathcal{M}'$.
\end{theorem}

Here, we chose the wording similarly as in \cite{BK_contractibility}, which seems sufficient for applications, but we remark that proof also yields polynomial dependence on the curvature scale similarly as in \cite[Theorem 1.7]{BK_uniqueness}.\\

To establish the above results we observe that most steps of the proof from Bamler-Kleiner \cite{BK_uniqueness} go through in arbitrary dimension, with the only exception of the steps that use Perelman's description of 3-dimensional $\kappa$-solutions \cite{Per2}, Brendle's uniqueness result for 3-dimensional noncollapsed steady Ricci solitons \cite{Brendle_steady}, or the Anderson-Chow estimate for 3-dimensional linearized Ricci-DeTurck flow \cite{AC}. We then argue how recent results by Brendle-Naff \cite{BN}, Brendle-Daskalopoulos-Naff-Sesum \cite{BDNS}, and 
Chen-Wu \cite{CW} can be used as appropriate substitutes for these important ingredients in higher dimensions.\\

Finally, motivated by our recent proof of the mean-convex neighborhood conjecture for mean curvature flow through neck-singularities \cite{CHH,CHHW} it seems reasonable to conjecture that the canonical neighborhood property actually follows from infinitesimal properties (as opposed to global curvature assumptions such as PIC), specifically that any complete Ricci flow spacetime with compact initial condition for whose metric completion all tangent flows at singular points, c.f. \cite{Bam2,Bam3}, are round shrinking spherical space forms  or round shrinking necks, i.e. $ \mathbb{R}\times S^{n-1}$, satisfies the canonical neighborhood property.\\


\noindent\textbf{Acknowledgments.} My research has been supported by an NSERC Discovery Grant (RGPIN-2016-04331) and a Sloan Research Fellowship, and I thank Alexander Kupers and Keaton Naff for useful discussions.

\section{The proofs}

We will now describe the necessary modifications to the proof from \cite{BK_uniqueness} to establish the stability of singular Ricci flows in higher dimensions:

\begin{proof}[{Proof of Theorem \ref{thm_stability}}]
In the proof of Bamler-Kleiner \cite{BK_uniqueness}, both the comparison domain $\mathcal{U}$ and the comparison map $\phi: \mathcal{U}\to \phi(\mathcal{U})$ are constructed via induction on time. To this end, given $\mathcal{M}$ and $\eps$, they fix a small enough comparison scale $r_{\textrm{comp}}>0$, and consider discrete time steps $t_j = jr_{\textrm{comp}}^2$. Each time-slab $\mathcal{U}^j:= \mathcal{U}\cap \{ t_{j-1}\leq\mathfrak{t} < t_{j}\}$ is chosen to simply be a product domain, and the evolution on $\mathcal{U}^j$ can be described by an ordinary Ricci flow parametrized for $t\in [t_{j-1},t_j)$. For each time $t_j$, they construct the initial time slice of $\mathcal{U}^j$ from the final time slice of $\mathcal{U}^{j-1}$, by removing certain parts via cutting along central spheres of some necks and/or by adding certain parts via gluing in Bryant caps. This adding and removing is done in a careful way to ensure that certain a priori assumptions are preserved. In particular, the initial time-slice of $\mathcal{U}^j$ contains all points with curvature scale at least $\Lambda r_{\textrm{comp}}$, and on the other hand all points in $\mathcal{U}^j$ have curvature scale bigger than $\lambda r_{\textrm{comp}}$, where $0<\lambda<\Lambda<\infty$ are suitable constants. The comparison map $\phi$ is constructed as a suitable solution of the harmonic map heat flow for its inverse $\psi:=\phi^{-1}$,
\begin{equation}\label{harm_map}
\partial_t \psi =\Delta_{g_t,g_t'} \psi.
\end{equation}
This approach of course goes back to DeTurck \cite{DT}, who observed that if $g$ and $g'$ evolve by Ricci flow, then the difference tensor
\begin{equation}
h(t):=\phi(t)^\ast  g'_t - g_t,
\end{equation}
which quantifies how much $\phi(t)$ deviates from being an isometry, satisfies an equation of the form
\begin{equation}
D_t h = \Delta h +2\mathrm{Rm}(h) + \nabla h\ast\nabla h + h\ast\nabla^2 h.
\end{equation}
To show that $h$ remains small, Bamler-Kleiner consider the quantity
\begin{equation}
Q:=e^{H(\eps^{-1}-\mathfrak{t})} \rho_1^E |h|,
\end{equation}
where $H$ and $E$ are large enough constants, and $\rho_1:=\min\{\rho,1\}$,  where $\rho$ denotes the curvature scale. Their key analytic estimate to control $h$ (and thus to ensure that $\phi$ remains an almost isometry) is the semi-local maximum principle \cite[Proposition 9.1]{BK_uniqueness}, which says that if the evolution around $x\in\mathcal{M}$ is described by ordinary Ricci flow in a sufficiently large parabolic ball $P:=P(x,L\rho_1(x))$, and if $\sup_{P}|h|\ll 1$, then
\begin{equation}\label{semi_loc}
Q(x)\leq \frac{1}{100}\sup_P Q.
\end{equation}
This estimate is crucial to justify that comparison domain $\mathcal{U}$ and the comparison map $\phi$  indeed satisfy all the a priori assumptions (APA1)--(APA13) listed in \cite[Section 7]{BK_uniqueness}.\\

Most of the proof of Bamler-Kleiner \cite{BK_uniqueness}, as reviewed above, goes through perfectly well in arbitrary dimensions (with only minor cosmetic changes such as replacing $S^2$ by $S^{n-1}$). The only two points that use the dimension in an essential way are the description of canonical neighborhoods and the proof of the semi-local maximum principle estimate \eqref{semi_loc}. Let us now describe how to address these two points:\\

In recent important work, Brendle-Naff \cite{BN} and Brendle-Daskalopoulos-Naff-Sesum \cite{BDNS} proved that any higher-dimensional $\kappa$-solution (as defined in the introduction) is up to finite quotients, scaling and isometries, either the round shrinking $S^n$, the round shrinking $\mathbb{R}\times S^{n-1}$, the $\mathrm{SO}(n)$-symmetric Bryant soliton, or the $\mathrm{SO}(n)$-symmetric Perelman oval. This is a more than sufficient substitute for the description of canonical neighborhoods via results from \cite{Per2,Brendle_steady} that have been used in \cite{BK_uniqueness}.\footnote{Coincidentally, at the time when the paper by Bamler-Kleiner \cite{BK_uniqueness} was written, even in dimension 3 the complete classification of $\kappa$-solutions, due to Brendle \cite{Brendle_ancient_Ricci} and Brendle-Daskalopoulos-Sesum \cite{BDS}, was not available, so even in dimension 3 some steps could be simplified a bit by using the now available complete classification of $\kappa$-solutions.}\\

Via a standard blowup argument \cite[Section 9]{BK_uniqueness}, which goes through in arbitrary dimensions, the proof of \eqref{semi_loc} can be reduced to establishing a Liouville theorem for ancient solutions of the linearized Ricci-DeTurck equation
\begin{equation}
D_t h = \Delta h +2\mathrm{Rm}(h)
\end{equation}
on $\kappa$-solutions. Specifically, it has been shown in \cite[Theorem 9.8]{BK_uniqueness} that for any $\chi>0$ one has
\begin{equation}\label{thm_li}
\sup_{M\times (-\infty,0] } \frac{|h|}{R^{1+\chi}} <\infty\qquad \Rightarrow \qquad h\equiv 0.
\end{equation}
This in turn followed from a result by Anderson-Chow \cite{AC}, which says that on any 3-dimensional Ricci flow with positive scalar curvature one has the evolution inequality
\begin{equation}\label{eq_AC}
\partial_t \left( \frac{|h|^2}{R^2}\right) \leq \Delta \left( \frac{|h|^2}{R^2}\right) + \frac{2}{R} \nabla R\cdot \nabla  \left( \frac{|h|^2}{R^2}\right),
\end{equation}
and consequently by the strong maximum principle if $|h|/R$ attain its maximum over $P=P(x,r)$ at $x$, then it is constant on $P$. The inequality \eqref{eq_AC} does not hold for general higher dimensional Ricci flows. However, by the above classifications all $\kappa$-solutions are rotationally symmetric, and thus have vanishing Weyl tensor. As observed by Chen-Wu \cite{CW}, the evolution inequality \eqref{eq_AC} holds for higher dimensional Ricci flows with vanishing Weyl tensor. In conclusion, the Liouville estimate \eqref{thm_li} also holds in higher dimensions.\\

Using the above modifications, the proof from Bamler-Kleiner \cite{BK_uniqueness} goes through in all dimensions.
\end{proof}

To conclude, let us briefly explain how the other results stated in the introduction follow:

\begin{proof}[{Proof of Theorem \ref{thm_uniq}}]
The stated uniqueness result now easily follows by applying Theorem \ref{thm_stability} (stability of singular Ricci flows) for $\eps_i=1/i$ and passing to a limit of the $\eps_i$-isometries.
\end{proof}

\begin{proof}[{Proof of Corollary \ref{cor_well}}]
Existence follows from Hamilton \cite{Hamilton_pic} and Chen-Zhu \cite{CZ}, together with \cite[Appendix A]{KL}, and uniqueness has been established in Theorem \ref{thm_uniq} (uniqueness of singular Ricci flows).
\end{proof}

\bibliography{uniqueness_sing_ricci.bib}
\bibliographystyle{alpha}

\vspace{10mm}

{\sc Robert Haslhofer, Department of Mathematics, University of Toronto,  40 St George Street, Toronto, ON M5S 2E4, Canada}\\

\end{document}